\documentclass[12pt]{article}

\usepackage{amsmath,amsthm}
\usepackage{amssymb}
\usepackage{color}
\usepackage{xspace}
\usepackage[colorlinks=true,
linkcolor=green,
filecolor=brown,
citecolor=green]{hyperref}

\def\c{Ces\`{a}ro}
\def\th{\theta}
\def\LHS{left hand side }

\newcommand{\StirlingPartition}[2]{\genfrac{ \{ }{ \} }{0pt}{}{#1}{#2}}

\begin{document}
\newtheorem{lemma}{Lemma}
\newtheorem{theorem}{Theorem}
\newtheorem{prop}{Proposition}
\newtheorem{cor}{Corollary}
\vspace*{-25mm}
\begin{center}
{\Large
Ces\`{a}ro's Integral Formula for the Bell Numbers (Corrected)                           \\ 
}

\vspace{10mm}
DAVID CALLAN  \\
Department of Statistics  \\
\vspace*{-2mm}
University of Wisconsin-Madison  \\
\vspace*{-2mm}
Medical Science Center \\
\vspace*{-2mm}
1300 University Ave  \\
\vspace*{-2mm}
Madison, WI \ 53706-1532  \\
{\bf callan@stat.wisc.edu}  \\
\vspace{5mm}

October 3, 2005
\end{center}

\vspace{5mm}

In 1885, Ces\`{a}ro \cite{cesaro} gave the remarkable formula
\[
N_{p}=\frac{2}{\pi e}\int_{0}^{\pi}e^{e^{\cos \theta}\cos(\sin 
\theta))} \sin(\,e^{\cos \theta} \sin(\sin \theta)\,) \sin p \theta\ 
d\theta
\]
where $(N_{p})_{p\ge 1}=(1,2,5,15,52,203,\ldots)$ are the modern-day 
Bell numbers. This formula  
was reproduced 
verbatim in the Editorial Comment on a 1941 Monthly problem 
\cite{becker} (the notation $N_{p}$ for Bell number was still in use 
then). I have not seen it in recent works and, while it's not very 
profound, I think it deserves to be 
better known.

Unfortunately, it contains a typographical error: a factor of $p!$ is 
omitted. The correct formula, with $n$ in place of $p$ and using $B_{n}$ for Bell number, is
\[
\hspace*{30mm}    B_{n}=\frac{2\,n!}{\pi e}\int_{0}^{\pi}e^{e^{\cos \theta}\cos(\sin 
\theta))} \sin(\,e^{\cos \theta} \sin(\sin \theta)\,) \sin n \theta\ 
d\theta \hspace*{10mm} n\ge 1.
\]
The integrand is the imaginary part of $e^{e^{e^{i\theta}}}\!\sin n\th$, 
and so an equivalent formula is
\begin{equation}
    B_{n}  =\frac{2\,n!}{\pi e} \textrm{Im} \left(
    \int_{0}^{\pi}e^{e^{e^{i\theta}}}\sin n\theta\ d\th \right).
    \label{eq:2}
\end{equation}
The formula (\ref{eq:2}) is quite simple to prove modulo a few 
standard facts about set partitions. Recall that the Stirling partition 
number $\StirlingPartition{n}{k}$ is the number of partitions of 
$[n]=\{1,2,\ldots,n\}$ 
into $k$ nonempty blocks and the Bell number 
$B_{n}=\sum_{k = 1}^{n}\StirlingPartition{n}{k}$ counts all partitions 
of $[n]$. Thus $k! \StirlingPartition{n}{k}$ counts ordered partitions 
of $[n]$ into $k$ blocks (the $k!$ factor serves to order the blocks) 
or, equivalently,  
counts surjective 
functions $f$ from $[n]$ onto $[k]$ (the $j$th block is $f^{-1}(j)$). Since the number of unrestricted 
functions from $[n]$ to $[j]$ is $j^{n}$, a classic application of the 
inclusion-exclusion principle yields
\begin{equation}
     k! \StirlingPartition{n}{k} = 
     \sum_{j=0}^{k}(-1)^{k-j}\binom{k}{j}j^{n}.
    \label{eq:3}
\end{equation}

The trig identity underlying \c's formula is nothing more than the 
orthogonality of sines on $[0,\pi]$:
\[
\int_{0}^{\pi}\sin m\theta \sin n\th\ d\th  =
\begin{cases}
    \frac{\mbox{\raisebox{.2ex}{$\pi$}}}{\mbox{\raisebox{-.5ex}{2}}} & \textrm{if }m=n, \\
    0 &  \textrm{if }m\ne n 
\end{cases}
\]
for $m,n$ nonnegative integers. Using the Taylor expansion $e^{x}=\sum_{m\ge 
0}\frac{x^{m}}{m!}$ and DeMoivre's formula $e^{i\th}=\cos\th+i \sin\th$, it follows that
\begin{equation}
    \textrm{Im}\left(\int_{0}^{\pi}e^{j e^{i\th}}\sin n\th \ 
    d\th\right)=\frac{j^{n}}{n!}\frac{\pi}{2}
    \label{eq:4}
\end{equation}
for integer $j\ge 0$. Just one more identity is needed:
\begin{equation}
  \textrm{Im}\left(\int_{0}^{\pi}\frac{\big(e^{e^{i\th}}-1\big)^{k}}{k!}
  \sin n\th\ d\th\right) =\frac{1}{n!}\StirlingPartition{n}{k}\frac{\pi}{2}
    \label{eq:5}
\end{equation}
for integer $k\ge 0$ (of course, $\StirlingPartition{n}{0}=0$ for 
$n>0$ and $\StirlingPartition{n}{k}=0$ for $k>n \ge 0$).

\textbf{Proof of (\ref{eq:5}) }\quad The binomial theorem implies the \LHS is
\begin{eqnarray*}
     &  & \frac{1}{k!} \sum_{j=0}^{k}(-1)^{k-j}\binom{k}{j}\textrm{Im}\left(\int_{0}^{\pi}
     e^{je^{i\th}}\sin n\th \ d\th\right)
      \\
     & \underset{(\ref{eq:4})}{=} & 
    \frac{1}{k!} \sum_{j=0}^{k}(-1)^{k-j}\binom{k}{j}\frac{j^{n}}{n!}\frac{\pi}{2}
     \\
     &  \underset{(\ref{eq:3})}{=} & 
     \frac{1}{n!}\StirlingPartition{n}{k} \frac{\pi}{2}
\end{eqnarray*} \qed

Summing (\ref{eq:5}) over $k \ge 0$ yields \c's formula 
(\ref{eq:2}). The Bell numbers have many other pretty representations, including 
Dobinski's infinite sum formula \cite[p.\,210]{comtet}
\[
B_{n}=\frac{1}{e}\sum_{k=0}^{\infty}\frac{k^{n}}{k!}.
\]

AMS Classification numbers: 05A19, 05A15.
\end{document}